\input amstex
\documentstyle{amsppt}
\magnification=\magstep1
 \hsize 13cm \vsize 18.35cm \pageno=1
\loadbold \loadmsam
    \loadmsbm
    \UseAMSsymbols
\topmatter
\NoRunningHeads
\title Note on multiple $q$-zeta functions
\endtitle
\author
  T. Kim
\endauthor
 \keywords : multiple $q$-zeta function, q-Euler numbers,
  Laurent series, Cauchy integral
\endkeywords

\abstract  In this paper we consider the analytic continuation of the multiple Euler $q$-zeta function
 in the complex number field as follows:
$$\zeta_{r,q}^E(s,x)=[2]_q^r\sum_{m_1,\cdots, m_r=0}^{\infty}\frac{(-1)^{m_1+\cdots+m_r}}{[x+m_1+\cdots+m_r]_q^s},$$
where $q\in \Bbb C$ with $|q|<1$, $\Re(x)>0$, and $r\in \Bbb N$.
Thus, we investigate their behavior near the poles and give the corresponding functional equations.
\endabstract
\thanks  2000 AMS Subject Classification: 11B68, 11S80
%\newline keywords and phrases
\newline  The present Research has been conducted by the research
Grant of Kwangwoon University in 2010
\endthanks
\endtopmatter

\document

{\bf\centerline {\S 1. Introduction/ Preliminaries}}

 \vskip 15pt
  Let $\Bbb C$ be the complex number field. For $s\in \Bbb C$, the Hurwitz's type Euler zeta function is defined
by
$$\zeta^E (s, x)=2\sum_{k=0}^{\infty}\frac{(-1)^k}{(k+z)^s}, \text{ where $s\in \Bbb C$, $z\neq 0, -1, -2, \cdots,$ ( see [11] ).}\tag1$$
  Thus, we note that $\zeta^{E}(s, x)$ is a meromorphic function in whole complex $s$-plane.  It is well known that the Euler polynomials are defined as
  $$\frac{2}{e^t +1}e^{xt}=\sum_{n=0}^{\infty}E_n(x)\frac{t^n}{n!}, \text{  for $|t|<\pi$,}\tag2$$
      and $E_n=E_{n}(0)$ are called the $n$-th Euler numbers (see [ 7, 8, 9, 11]).  By (1) and (2), we note that
      $\zeta^E(-n, x)=E_n(x)$, for $n\in \Bbb Z_+$.  Throughout this paper we assume that $q\in \Bbb C$ with $|q|<1$ and we use
the notation of $q$-numbers as $ [x]_q=\frac{1-q^x}{1-q}.$
    The $q$-Euler numbers are defined as
    $$E_{0,q}=\frac{2}{[2]_q}, \text { and } (qE+1)^n+E_{n,q}=0 \text{ if $n\geq 1$ }  ,\tag3 $$
    where we use the standard convention about replacing $E^k$ by $E_{k,q}$( see [7]). Thus, we define the $q$-Euler polynomials as follows:
$$ E_{n,q}(x)=\sum_{l=0}^n\binom{n}{l}q^{lx}[x]_q^{n-l}q^{lx}E_{l,q}, \text{ (see [7, 8, 15])}.\tag4$$
For $s\in \Bbb C$, the q-extension of Hurwitz's type $q$-Euler zeta function is defined by
$$\zeta_q^E (s,x)=[2]_q\sum_{n=0}^{\infty}\frac{(-1)^n}{[n+x]_q^s}, \text{ where $x\neq 0, -1, -2,\cdots$.}\tag5$$
For $n\in\Bbb Z_{+}$, we have $\zeta_q^E(-n,x)=E_{n,q}(x)$(see [6, 7, 15] ). Let $\chi$ be a Dirichlet's character with conductor $f\in\Bbb N$ with $f\equiv 1$$( mod  \  2 )$. It is known that the generalized $q$-Euler polynomials attached to $\chi$ are defined by
$$F_{q,\chi}(t,x)=[2]_q\sum_{m=0}^{\infty}(-1)^m\chi(m)e^{[m+x]_qt}=\sum_{m=0}^{\infty}E_{m,\chi,q}(x)\frac{t^m}{m!}, \text{ see [7] }.\tag6$$
Note that $$\lim_{q\rightarrow 1}F_{q,\chi}(t,x)=\frac{ 2\sum_{a=1}^{f-1}(-1)^a \chi(a)e^{at}}{e^{ft}+1} e^{xt}=\sum_{m=0}^{\infty}E_{m,\chi}(x)\frac{t^m}{m!},$$
where $E_{m,\chi}(x)$ are called the $m$-th generalized Euler polynomials attached to $\chi$.
 From (6), we can derive the following equation.
 $$E_{n,\chi,q}(x)=[2]_q\sum_{m=0}^{\infty}(-1)^m\chi(m)[m+x]_q^n
 =\frac{[2]_q}{(1-q)^n}\sum_{l=0}^n\binom{n}{l}(-q^x)^l\sum_{a=0}^{f-1}\frac{(-1)^a\chi(a)q^{la}}{1+q^{lf}}.\tag7$$
 Now, we consider the Dirichlet's type Euler $q$-$l$-function which interpolate $E_{n,\chi,q}(x)$ at negative integer.
For $s\in\Bbb C$, define
$$ l_q(s,x|\chi)=\sum_{n=0}^{\infty}\frac{\chi(n)(-1)^n}{[n+x]_q^s}, \text{ $x\neq 0, -1, -2,\cdots$, (see [6, 7, 8, 15]). } $$
Note that $l_q(-n, x|\chi)=E_{n,\chi,q}(x)$ for $n\in\Bbb Z_+$. In the special case $x=0$, $E_{n,\chi,q}(=E_{n,\chi,q}(0))$ are called
the $n$-th generalized Euler numbers attached to $\chi$.
The theory of quantum groups has been quite successful in producing identities for $q$-special function. Recently, several mathematicians
have studied $q$-theory in the several areas(see [1-23]). In this paper we approach the $q$-theory in the area of special function. That is,
we first consider the analytic continuation of multiple $q$-Euler zeta function in the complex plane as follows:
$$\zeta_{r,q}^E(s,x)=[2]_q^r\sum_{m_1,\cdots, m_r=0}^{\infty} \frac{(-1)^{m_1+\cdots+m_r}}{[x+m_1+\cdots+m_r]_q^s},
\text{ $s\in\Bbb C$, $x\neq 0, -1, \cdots.$}\tag8$$
From (8), we investigate some identities for the multiple $q$-Euler numbers and polynomials. Finally, we give interesting functional equation related to the multiple $q$-Euler polynomials, gamma functions and multiple $q$-Euler zeta function.

\vskip 10pt

{\bf\centerline {\S 2. Multiple $q$-Euler polynomials and multiple $q$-Euler zeta functions}} \vskip 10pt

 From (3), we note that
 $$E_{n,q}=\frac{[2]_q}{(1-q)^n}\sum_{l=0}^n\binom{n}{l}\frac{(-q^x)^l}{(1+q^l}
 =[2]_q\sum_{m=0}^{\infty}(-1)^m[m+x]_q^n.\tag9$$
 Let $F_q(t,x)=\sum_{n=0}^{\infty}E_{n,q}(x)\frac{t^n}{n!}.$ Then we see that
 $$F_q(t,x)=[2]_q\sum_{m=0}^{\infty}(-1)^me^{[m+x]_qt}. \tag10$$
From (10), we note that $\lim_{q\rightarrow 1}F_q(t,x)=\frac{2}{e^t+1}e^{xt}=\sum_{n=0}^{\infty}E_n(x)\frac{t^n}{n!},$  where $E_n(x)$
are called the $n$-th Euler polynomials. For $s\in \Bbb C$, we have
$$\frac{1}{\Gamma(s)}\int_{0}^{\infty}t^{s-1}F_q(-t,x)dt=[2]_q\sum_{n=0}^{\infty}\frac{(-1)^n}{[n+x]_q^s},
\text{ where $x\neq 0, -1, -2, \cdots.$}\tag11$$ By Cauchy residue theorem and Laurent series, we see that
$\zeta_q^E(-n, x)=E_{n,q}(x)$ for $n\in\Bbb Z_{+}.$
 Let $\chi$ be the Dirichlet's character with conductor $f(=odd)\in \Bbb N$. From (6), we can derive
 $$\aligned
 F_{q,\chi}(t,x)&=[2]_q\sum_{a=0}^{f-1}(-1)^a\chi(a)\sum_{n=0}^{\infty}(-1)^ne^{[a+x+nf]_qt}\\
 &=[2]_q\sum_{a=0}^{f-1}(-1)^a\chi(a)\sum_{n=0}^{\infty}(-1)^ne^{[f]_q[\frac{x+a}{f}+n]_{q^f}t}. \endaligned\tag12$$                                              Let us define the Dirichlet's type  $q$-Euler  $l$-function as follows:
$$l_q(s, x|\chi)=[2]_q\sum_{m=0}^{\infty}\frac{\chi(m)(-1)^m}{[m+x]_q^s}, \text{ where $s\in\Bbb C$, $x\neq 0, -1, -2,\cdots.$ }\tag13$$
From the Mellin transformation of $F_{q,\chi}(t,x)$, we note that
$$\frac{1}{\Gamma(s)} \int_{0}^{\infty}F_{q,\chi}(-t,x)t^{s-1}dt=[2]_q\sum_{n=0}^{\infty}\frac{(-1)^n\chi(n)}{[n+x]_q^s},
 \text{ where $s\in\Bbb C$, $x\neq 0, -1, -2,\cdots.$}\tag14$$
By Laurent series and Cauchy residue theorem, we see that
$ l_q(-n,x|\chi)=E_{n,\chi,q}(x)$ for $n\in \Bbb Z_+$. Let us consider the following $q$-Euler polynomials of order $r(\in \Bbb N$).
$$F_q^{(r)}(t,x)=[2]_q^r\sum_{m_1,\cdots, m_r=0}^{\infty}(-1)^{m_1+\cdots+m_r}e^{[m_1+\cdots+m_r+x]_qt}
=\sum_{n=0}^{\infty}E_{n,q}^{(r)}(x)
\frac{t^n}{n!}. \tag15$$
In the special case $x=0$, $E_{n,q}^{(r)}(=E_{n,q}^{(r)}(0))$ are called the $n$-th $q$-Euler numbers of order $r$.
It is easy to show that $\lim_{q\rightarrow 1}F_q^{(r)}(t,x)=\left(\frac{2}{e^t +1}\right)^r e^{xt}=\sum_{n=0}^{\infty}E_n^{(r)}(x)\frac{t^n}{n!},$
where $E_n^{(r)}(x)$ are called the $n$-th Euler polynomials of order $r$.
From (15), we note that
$$\aligned
 \sum_{n=0}^{\infty}E_{n,q}^{(r)}(x)\frac{t^n}{n!}&=[2]_q^r\sum_{m_1,\cdots, m_r=0}^{\infty}(-1)^{m_1+\cdots+m_r}e^{[m_1+\cdots+m_r+x]_qt}\\
 &=[2]_q^r \sum_{m=0}^{\infty}\binom{m+r-1}{m}(-1)^m e^{[m+x]_qt}.
 \endaligned\tag16$$
 Thus, we have
 $$E_{n,q}^{(r)}(x)=\frac{[2]_q^r}{(1-q)^n}\sum_{l=0}^n\binom{n}{l}(-1)^lq^{lx}\left(\frac{1}{1+q^l}\right)^r.$$
 Therefore, we obtain the following proposition.

\proclaim{ Proposition 1} For $r\in\Bbb N$,  $n\in \Bbb Z_+$, we have
$$\aligned
 E_{n,q}^{(r)}(x)&=[2]_q^r\sum_{m_1, \cdots,
m_r=0}^{\infty}(-1)^{m_1+\cdots+m_r}[m_1+\cdots+m_r+x]_q^n\\
&=[2]_q^r\sum_{m=0}^{\infty}\binom{m+r-1}{m}(-1)^m [m+x]_q^n\\
&=\frac{[2]_q^r}{(1-q)^n}\sum_{l=0}^n\binom{n}{l}(-1)^lq^{lx}\left(\frac{1}{1+q^l}\right)^r.
\endaligned$$
\endproclaim
 By Mellin transformation of $F_q^{(r)}(t,x)$, we see that
 $$\aligned
& \frac{1}{\Gamma(s)}\int_{0}^{\infty}F_q^{(r)}(-t,x)t^{s-1}dt=[2]_q^r\sum_{m=0}^{\infty}
 \frac{\binom{m+r-1}{m}(-1)^m}{[m+x]_q^s}\\
 &=[2]_q^r\sum_{m_1,\cdots, m_r=0}^{\infty}\frac{(-1)^{m_1+\cdots+m_r}}{[m_1+\cdots+m_r+x]_q^s}, \text{ where $s\in\Bbb C$, $x\neq 0,-1, -2, \cdots.$}
 \endaligned\tag17$$
From (17), we can consider the following multiple $q$-Euler zeta function.
\proclaim{ Definition 2}
 For $s\in\Bbb C$, $x\in\Bbb R$ with $x\neq 0,-1,-2, \cdots, $ we define the multiple $q$-Euler zeta function as follows:
 $$\zeta_{r,q}^E(s,x)=[2]_q^r\sum_{m_1,\cdots, m_r=0}^{\infty}\frac{(-1)^{m_1+\cdots+m_r}}{[m_1+\cdots+m_r+x]_q^s}.$$
\endproclaim
Note that $\zeta_{r,q}^E$ is meromorphic function in whole complex $s$-plane. By using Cauchy residue theorem and Laurent series in (15) and (17),
we obtain the following theorem.
\proclaim{ Theorem 3}
For $n\in\Bbb Z_+$, $r\in\Bbb N$, we have
$$\zeta_{r,q}^E(-n,x)=E_{n,q}^{(r)}(x).$$
\endproclaim
In (15), we have
$$F_q^{(r)}(t,x)=[2]_q^r\sum_{a_1,\cdots, a_r=0}^{f-1}(-1)^{a_1+\cdots+a_r}\sum_{m_1,\cdots, m_r=0}^{\infty}(-1)^{m_1+\cdots+m_r}e^{[\sum_{i=1}^r(a_i+fm_i)+x]_qt}. \tag18$$
By (17) and (18), we obtain the following theorem.
\proclaim{ Theorem 4} ( Distribution relation for $E_{m,q}^{(r)}(x)$)

For $n\in\Bbb Z_+$, $f\in\Bbb N$ with $f\equiv 1$ $(mod  \  2)$, we have
$$E_{n,q}^{(r)}(x)
=\left(\frac{[2]_q}{[2]_{q^f}}\right)^r [f]_q^n\sum_{a_1,\cdots, a_r=0}^{f-1}(-1)^{a_1+\cdots+a_r}
E_{n, q^f}^{(r)}\left(\frac{a_1+\cdots+a_r+x}{f}\right).$$
Moreover,
$$E_{n,q}^{(r)}(x)=\frac{[2]_q^r}{(1-q)^n}\sum_{l=0}^n\binom{n}{l}(-q^x)^l
\sum_{a_1,\cdots, a_r=0}^{f-1}\frac{(-1)^{a_1+\cdots+a_r}q^{l(a_1+\cdots+a_r)}}{(1+q^{lf})^r}.$$
\endproclaim
Let $\chi$ be the Dirichlet's character with conductor $f(=odd)\in \Bbb N$. Then we define the generalized $q$-Euler polynomials of
order $r$ attached to $\chi$ as follows:
$$\aligned
&F_{q,\chi}^{(r)}(t,x)=\sum_{n=0}^{\infty}E_{n,\chi,q}^{(r)}(x)\frac{t^n}{n!}\\
&=[2]_q^r\sum_{m_1,\cdots, m_r=0}^{\infty}(-1)^{ m_1+\cdots+m_r}\left( \prod_{j=1}^r\chi(m_j) \right)
e^{[x+m_1+\cdots+m_r ]_qt}.
\endaligned \tag19$$
In the special case $x=0$, $E_{n,\chi,q}^{(r)}(=E_{n,\chi,q}^{(r)}(0))$ are called the $n$-th generalized $q$-Euler numbers of order $r$ attached to $\chi$. From (19), we can derive
$$ \aligned
&F_{q,\chi}^{(r)}(t,x)=[2]_q^r\sum_{m_1,\cdots, m_r=0}^{\infty}(-1)^{ m_1+\cdots+m_r}\left( \prod_{j=1}^r\chi(m_j) \right)
e^{[x+m_1+\cdots+m_r ]_qt}\\
&=[2]_q^r\sum_{m=0}^{\infty}\binom{m+r-1}{m}(-1)^m\sum_{a_1,\cdots, a_r=0}^{f-1}\left(\prod_{j=1}^r \chi(a_j)\right)
(-1)^{\sum_{j=1}^r a_j}e^{[x+mf+\sum_{j=1}^ra_j]_qt}.\endaligned\tag20$$
By (16) and (20), we obtain the following theorem.
 \proclaim{ Theorem 5} For $f\in \Bbb N$ with $f\equiv 1$$(mod   \  2 )$, we have
$$ E_{n,\chi,q}^{(r)}(x)=[f]_q^n\left(\frac{[2]_q}{[2]_{q^f}}\right)^r \sum_{a_1,\cdots, a_r=0}^{f-1}\left(\prod_{j=1}^r\chi(a_j)\right)(-1)^{\sum_{j=1}^ra_j}E_{n,q^f}^{(r)}(\frac{x+\sum_{j=1}^ra_j}{f}),$$
and
$$E_{n,\chi,q}^{(r)}(x)=\frac{[2]_q^r}{(1-q)^n} \sum_{l=0}^n\binom{n}{l}(-q^x)^l\sum_{a_1,\cdots, a_r=0}^{f-1}
\frac{\left(\prod_{j=1}^r\chi(a_j)\right)(-q^l)^{\sum_{i=1}^ra_i}}{(1+q^{lf})^r}.$$
  \endproclaim
 From the Mellin transformation of $F_{q,\chi}^{(r)}(t,x)$, we note that
 $$\frac{1}{\Gamma(s)}\int_{0}^{\infty}F_{q,\chi}^{(r)}(-t,x)t^{s-1}dt
 =[2]_q^r\sum_{m_1,\cdots, m_r=0}^{\infty}\frac{(-1)^{m_1+\cdots+m_r}\left(\prod_{j=1}^r \chi(m_j))\right)}{[m_1+\cdots+m_r+x]_q^s},  \tag21$$                                   where $ s\in \Bbb C$, $ \Re(x) >0$. From (21) we can also consider the following Dirichlet's type multiple $q$-Euler $l$-function.
\proclaim{ Definition 6} For $s\in \Bbb C$, $x\in\Bbb R$ with $x\neq 0, -1, -2, \cdots,$
we define Dirichlet's type $q$-Euler $l$-function as follows:
$$l_q^{(r)}(s, x|\chi)=[2]_q^r\sum_{m_1,\cdots, m_r=0}^{\infty}\frac{(-1)^{m_1+\cdots+m_r} \left(\prod_{j=1}^r\chi(m_j)\right)}
{[m_1+\cdots+m_r+x]_q^s}.$$
\endproclaim
Note that $l_{q}^{(r)}(s,x|\chi)$ is also holomorphic function in whole complex $s$-plane. By (20) and (21), we see that
$$\aligned
&l_q^{(r)}(s, x|\chi)\\
&=\frac{1}{[f]_q^s}\left(\frac{[2]_q}{[2]_{q^f}}\right)^r\sum_{a_1,\cdots, a_r=0}^{f-1}
\left(\prod_{j=1}^r \chi(a_j)\right)(-1)^{\sum_{i=1}^ra_i}\zeta_{r,q^f}^E(s, \frac{a_1+\cdots+a_r+x}{f}).
\endaligned$$
By using Laurent series and Cauchy residue theorem, we obtain the following theorem.
\proclaim{ Theorem 7}
   For $n\in\Bbb Z_+$, we have
   $$l_q^{(r)}(-n, x|\chi)=E_{n,\chi, q}^{(r)}(x).$$
\endproclaim
For $q=1$, Theorem 7 seems to be similar type  of Dirichlet's $L$-function in complex analysis. That is, let $\chi$ be the
Dirichlet's character with conductor $d\in \Bbb N$. Then Dirichlet $L$-function is defined as
$$L(s,x|\chi)=\sum_{n=0}^{\infty}\frac{\chi(n)}{(n+x)^s}, \text{ where $s\in \Bbb C$, $x\neq 0, -1, -2, \cdots.$}$$
Let $n$ be positive integer. Then we have $L(-n,x|\chi)=-\frac{B_{n,\chi}(x)}{n},$ where $B_{n,\chi}(x)$ are called
the $n$-th generalized Bernoulli polynomials attached to $\chi$ (see [13, 14, 16, 18, 2, 3, 20-23]).

\vskip 20pt

\quad Taekyun Kim

\quad Division of General Education-Mathematics,

\quad Kwangwoon University,

\quad Seoul 139-701, S. Korea \quad e-mail:\text{
tkkim$\@$kw.ac.kr}

\enddocument